\documentclass[12pt]{article}
\usepackage[namelimits,intlimits]{amsmath}
\usepackage[english,russian]{babel}
\usepackage{amsfonts}
\usepackage{euscript}
\usepackage{amsthm}
\usepackage{multicol,amsthm,syntonly}
\usepackage{euscript}
\usepackage{color}
\usepackage{graphicx}
\usepackage{eufrak}
\begin{document}
\begin{center}

\vspace{30mm} {\Large \bf Linear quasigroups. II}\\
\vspace{10mm} {\bf Tabarov Abdullo\footnote[1]{The research to
this article was sponsored by Special Projects Office, Special
and Extension Programs of the Central European University Corporation. Grant Application Number: 321/20090.}}\\
{\it Tajik National University, Department of Mechanics and
Mathematics,\\
734017, Dushanbe, Rudaki ave.17, Tajikistan}\\
e-mail: tabarov2010@gmail.com
\end{center}

\vspace{20mm} \pagestyle{plain}
\setcounter{page}{1} \sloppy
\begin{center}
{\bf Abstract}
\end{center}

The article is a continuation of the author's work "Linear
quasigroups. I"  \ and devoted to linear quasigroups and some of
their generalizations. In the  second  part  identities  and
linearity of  quasigroups are investigated, in particular,  the
approach of finding quasigroup identities from some variety of
loops is  presented. Besides that we give  short chronology of
development the theory of quasigroups and loops and some
historical aspects.
\\



\vspace{20mm}

\begin{center}
{\bf I. Identities  and a linearity of  quasigroups}
\end{center}

First note that all necessary  definitions and  notions can be
found in [1, 2].

It is well known that the   theory of identities in  various
algebraic systems plays an  important role. In the theory of the
quasigroups many classes are given by identities (for example,
medial quasigroups, distributive quasigroups, Steiner quasigroups,
CH-quasigroups,  F-quasigroups etc.).

The theory of identities in the algebras has two interrelated
aspects: \emph{identities and algebras.} Accordingly, there are
two global questions (see [3]):

(1) describe  algebras with identities;

(2) describe  identities in algebras.

In the first case we have  a structure of algebra, in which  an
identity or systems of  identities of a special form is satisfied.
Often, on the structure of algebras noticeably  affected only
execution of a nontrivial identity, for example  in  the situation
of the associative rings. There are excellent  examples when  the
algebra of a certain natural class is completely  determined by
its identities.

The second case  is  about  finding the specific identities of
algebras or a class  of algebras. It should be noted that in this
formulation the  problem is  extremely difficult. For example, in
many cases there is no clarity in   question  of  existence of the
finite basis of a system of identities (\emph{the problem of
finite basis}). Moreover, the problem of identities, as a rule, is
separated from themselves algebras, in which these identities are
satisfied, although, of course, for each system of identities  a
class of algebras can be  chosen for which this system is defined.
Here we will discuss  varieties of algebras. The language of
varieties, in fact erases  a difference between algebras and
identities, or rather, it  allows easily  move from one of these
concepts to another. Precisely in  this language, as a rule, the
work in the theory of identities is provided.

To this we shall add that namely in  the second part of the
question G.B. Belyavskaya and A.Kh. Tabarov in [4-6] characterized
the classes of linear, alinear, mixed type of linearity, onesided
linear   and alinear quasigroups, $T$-quasigroups and some of
their generalizations by  identities and system of identities.

One of the properties  of the quasigroup which  characterizes its
nearness to  groups,  is availability of its  isotopy to  groups.
The class of quasigroups which are  isotopic  to  groups was
researched by many algebraists. This class is near to  groups and
so their construction  largely is easy for investigation. However,
in this class  there are many  difficult problems which are still
unresolved.

The class of quasigroups are isotopic to groups, first was
investigated by V.D. Belousov  [7]. In particular, V.D. Belousov
proved that the class of quasigroups isotopic to groups is
characterized by the following identity of  five variables:
$$
x\left( {y\backslash \left( {\left( {z / u\left. \right)v\left.
\right)} \right.} \right.} \right) = \left( {\left( {x\left(
{y\backslash z\left. {\left. \right)} \right) / u\left. \right)v}
\right.} \right.} \right..
$$
In  other words a variety  of quasigroups are isotopic to groups
characterized by this identity. Later, F.N.Sokhatsky noticed that
quasigroups which are isotopic to groups, can be characterized by
the identity with four variables [8]:
$$
[(x(u\backslash y)) / u]z = x[u\backslash ((y / x)z)].
$$

A variety are also all quasigroups isotopic to  abelian groups. It
also was first noted by V.D.Belousov [7]. This class of
quasigroups characterized by the identity of the four variables:
$$
x \backslash (y(u\backslash v))=u\backslash(y(x\backslash v)).
$$

In [2] we have given some information about  linear quasigroups
and their generalizations. Recall that a quasigroup $(Q,\cdot)$ is
called {\it linear} over a group $(Q, + )$, if it has the form:
\begin{equation}
\label{eq1} xy = \varphi x + c + \psi y,
\end{equation}
where $\varphi ,\psi \in Aut(Q, + )$,  $c$ is a fixed element from
$Q$.

By analogy with linear quasigroups,  G.B. Belyavskaya and A.KH.
Tabarov in [4,5]  defined alinear quasigroups, as well as
introduced the classes of left and right linear quasigroups, left
and right alinear quasigroups and mixed type of linearity.

In any algebra  identities are closely connected with algorithmic
problems. Articles of T. Evans [9,10], M.M. Glukhov and A.A.
Gvaramiya [11] were devoted to the algorithmic problems of the
theory of quasigroups. In 1951 T. Evans [9] proved an assertion
about positive solvability of word problem  for finitely-presented
algebras of all variety algebras $V(\Sigma )$, in which  the
theorem about embedding of such extremity partial $(\Sigma
)$-algebra in algebra from $V$ has a place. Further, by the way of
carried words to canonical type T. Evans  proved the isomorphism
problem for some classes of multiplicative system, in particular
for quasigroups and loops [10].

In 1970-1971  M.M. Glukhov [12] formulated the condition which is
stronger  then the theorem of embedding  and  such  that in
variety not only quasigroups, but in universal algebras  the
algorithmic problems of word, isomorphism and embedding are
positively  solved. M.M. Glukhov  called this condition $R$. The
varieties on algebra in which the condition $R$ is satisfied, were
called $R$-varieties.

In the research of algebras of one variety, sometimes it should be
helpful to use algebras of other variety  concerned with it.
Exactly this method was used in a series of  works during the
research of quasigroups isotopic to groups. In this connection it
is important to consider the approach know as equivalence and
rational equivalence of a class of algebras presented in [13]
which we have used in proving  the main theorem. Using  the
notions  of an equivalence and a rational equivalence in [14] was
proved that class primitive linear quasigroups, and in particular,
$T$-quasigroups are varieties. G.B. Belyavskaya and A.Kh. Tabarov
in [5,6] established  the various systems of identities which
characterized  some varieties of linear quasigroups and
T-quasigroups.

We present here the meaning of the equivalence class of algebras
introduced in [13].

{\bf Definition 1.} {\it The class of algebras $K_{1}, K_{2}$ with
accordingly to signatures $\Delta_1, \Delta_2$ are called
equivalent if the bijective  map
$$
f: K_{1} \to K_{2},
$$
exists, satisfying the conditions:

1) For any algebra $A \in K_{1}$ the  basis  algebras $A$ and
$f(A)$ coincide;

2) For any algebras $A, B\in K_{1}$ the map of $A$ in $B$ is a
homomorphism if and only if it is a homomorphism of algebra $f(A)$
in $f(B)$.}

At the same time the map $f$ is called  {\it an equivalence}
between classes $K_{1}$ and $K_{2}$.

If  $f$ is an  equivalence  between the varieties  of algebras
$K_{1}$ and $K_{2}$, then   the following assertions are hold:

{\it - Subset $A_{1} \subset A$  is   the subalgebra of algebra A,
if and only if $f(A_{1})$ is the subalgebra of algebra $f(A)$.

- Subset $M \subset A$   generates the algebra A if and only if
$M$ generates the algebra  $f(A)$.

- An equivalent relation in A is  a congruence of algebra A if and
only if it is a congruence of algebra  $f(A)$.

- Algebra $A \in K_{1}$  is free in variety  $K_{1}$ with the
basis $M$ if and only if $f(A)$ is free  in $K_{2}$ with  the
basis $M$.

- A class of algebras $L$  from  $K_{1}$  is  a subvariety in
$K_{1}$ if and only if the class  ${\{}f(A): A \in L{\}}$ is a
subvariety in $K_{2}$.}

Among all the classes of equivalences    of algebras, rational
equivalences particulary stand out.

Let  $K$  be the class of  algebras of  signature $\Delta$  and
$X$ be the alphabet of variables. Then on each $\Delta$-word $P$
in alphabet $X$ containing in its entries exactly $n$ variables,
for example $x_{1},\ldots ,x_{n}$,  the n-ary operation $w_{P}$ on
every algebra $A$ from $K$ could be determined. The meaning of
this operation on elements $a_{1},\ldots, a_{n} \in A$  is equal
to the meaning in $A$ of a word received from $P$ replacing
$x_{1},\ldots, x_{n}$ respectively with the elements
$a_{1},\ldots, a_{n}$. A such determined operation is called  {\it
a derivative operation} in signature $\Delta$ corresponding to the
$\Delta$-word $P$.

The transfer of the signature  $\Delta _{1}$ into signature
$\Delta_{2}$ is called any reflection $\tau$: $\Delta_{1} \to
\Delta_{2}$, such that every operation $h \in \Delta_{1}$ is
reflected in some derivative operation in the same arity  in
signature $\Delta_{2}$. If $\tau$ is such a transfer then on any
algebra $B$ of signature $\Delta_{2}$ the algebra $A=T_\tau (B)$
of signature $\Delta_{1}$ could be determined  with the same main
set if an  operation $h \in \Delta_1$ which coincides  with
$\tau$(h) could also be determined.

\emph{The class of algebras $K_{1}$, $K_{2}$  with the respective
signatures $\Delta_1, \Delta_2$ is  called rationally equivalent
if  the transfer of signatures $\tau$: $\Delta_{1} \to \Delta_{2}$
and $\sigma$: $\Delta_{2} \to \Delta_{1}$ and the bijective
mapping $f: K_{1} \to K_{2}$ exists such  that $T_\tau T_\sigma$
and $T_\sigma T_\tau$ are identity mappings.}

It is clear that a rational equivalence class of the algebras  is
an  equivalence and therefore for rational equivalence  classes
the above introduced  assertion is true.

In variety of quasigroups as well as for all algebras, a free
quasigroups of this variety plays the important role. However for
all variety of algebras to constructively describe free quasigroup
with the basis $A$ in this or that variety of quasigroup it is
necessary to have an algorithm of recognition of equivalence  of
words. In the general case this problem is very difficult. The
solution of this problem in given variety $Q(\Omega, \Sigma)$,
naturally depends on determination of its system of identities
$\Sigma$. There are varieties of quasigroups with solved problem
of words. All $R$-varieties [11,12] are related to these
varieties. At the same time there exists varieties with unsolved
problem of words in free quasigroups [15].

For any variety of quasigroups   there is an interesting question
about a connection of its free quasigroups with the free groups of
its isotopy closing  group. First  this problem for the variety of
T-quasigroups was considered in [16,17]. In these works  the
construction of free T-quasigroups, were offered based on the use
of the concept rational equivalence  of varieties of algebra with
different signatures were introduced in [14](see also [18]).

The class of quasigroups $Q(\Gamma)$ in the signature
$\Omega_{1}={\{}\cdot, /, \backslash, u{\}}$, where $u$ is  symbol
of $0$-ary operation and variety of $U$ algebras of signature
$\Omega_{2}={\{}+, -, 0,$ $\alpha, \beta, \gamma, \delta, c{\}},$
where $\alpha, \beta, \gamma, \delta$ symbols of unary and
\emph{c}  symbol of $0$-ary operation are given by the system of
group identities
$$
(x+y)+z=x+(y+z), x+0=x, 0+x=x, x+(-x)=0, (-x)+x=0
$$
in signature ${\{}+, -, 0{\}}$ and identities:
$$
x\alpha \gamma = x\gamma \alpha = x,\,\, x\beta \delta = x\delta
\beta = x,\,\, 0\alpha = 0\beta = 0\gamma = 0\delta = 0.
$$
was considered in [14].

Below we  substituted in the last system $\gamma$ into
$\alpha^{-1}$ and $\delta$ into $\beta^{-1}$, and will write it in
the next form:
$$
x\alpha \alpha^{-1} = x\alpha^{-1}\alpha = x,\,\, x\beta
\beta^{-1}=x\beta^{-1}\beta=x,\,\,
0\alpha=0\beta=0\alpha^{-1}=0\beta^{-1}=0.
$$
As transfers $\tau$: $\Omega_1 \to \Omega_2$ and $\sigma: \Omega_2
\to \Omega_1$  were taken mappings:
$$
\tau (\cdot)(x,y)=x\alpha + d + y\beta,
$$
$$
\tau (/)(x,y)=(x - y\beta - c)\gamma,
$$
$$
\tau(\backslash)(x,y) =(-c - x\alpha + y)\delta,
$$
$$
\tau(u)=0,
$$
$$
\sigma(+)(x,y)=(x/u)((u/u)\backslash y),
$$
$$
\sigma (-)(x)=(u/u)((x/u)\backslash u),
$$
$$
\sigma(\alpha )=x(u\backslash u),\,\,\, \sigma(\beta)=(u/u)x,
$$
$$
\sigma(\gamma) = x/(u\backslash u),\,\,\,
\sigma(\delta)=(u/u)\backslash x,\
$$
$$
\sigma (0)=u, \sigma (c)=uu.
$$

It was proved, that the classes of algebras $Q(\Gamma)$  and $U$
are rationally equivalence. From here it follows, in particular,
that $Q(\Gamma)$  are variety quasigroups. Analogously it was
proved that the varieties are the classes of right-linear,
left-linear, the linear quasigroups and $T$-quasigroups.

The authors observe that the received results about rational
equivalence allows   many questions about quasigroups  could  be
formulated in the  more usual language of algebras   from the
varieties $U$ which are near to groups. In particular it is noted
that in this way one can construct  free quasigroups from the
variety quasigroups $Q(\Gamma)$,$Q(A\Gamma)$ ($Q(A\Gamma)$ is the
variety of quasigroups isotopic to abelian groups). This idea,
with implicit use  of the equivalence  of varieties, was realized
in [16,17] for  $T$-quasigroups. They constructed a free algebra
with basis $X$ in some  different terms in the signature which
received a broader  signature by the unary operation symbols
$\alpha, \alpha^{-1}, \beta, \beta^{-1}$, on which were determined
the quasigroup operations and  in  result  a  free $T$-quasigroup
was received. Following attentive analysis of the works [16,17] it
became clear that some of the analogies of this construction
should be present for the variety of  all linear quasigroups [19].

\quad As it was noted above the constructive description of the
free quasigroups requires the decision  problem of word in the
free quasigroups in accordance of variety or as well as problem of
identities relations for the quasigroups of the considered
varieties. In this way  M.M. Glukhov in the report on the
algebraic conference devoted to 100-anniversary of A.G. Kurosh
announced the next result: if for  free groups some  varieties of
group solvable the word problem, then analogous fact is true for
its isotopic closure [20]. However from this result does not
follows the  decision of  word problem for  varieties of different
types of  linear quasigroups, as far as in these varieties
essential  role plays  limitations on the permutations which are
components of isotopy. Therefore the word problem for free
quasigroups from varieties of different types of linear
quasigroups  remains open.

Nevertheless  for some subclasses of  linear quasigroups, namely
$T$-quasigroups and medial quasigroups  the word problem for a
free  algebras is  solveable [21].

{\bf Theorem 1.1}. {\it In a variety of  all $T$-quasigroups  the
word problem   for  free algebras is solved}.

{\bf Proof:} To prove this  theorem we will need use the auxiliary
variety of algebras, namely the  variety $U(\Delta_1, S)$ algebras
of signature $\Delta_1 = \Delta \cup \Delta_0$   with system of
identities $S$, where $\Delta_0 = \{\alpha, \alpha^{-1}, \beta,
\beta^{-1}{\}}$ is the system of 0-ary operations, $\Delta ={\{}+,
-, 0{\}}$ is the group signature, and $S$ is the   system of
identities:
\begin{equation}
\begin{split}
\left. {\begin{array}{l} \label{4.5.1} (x+y)+z=x+(y+z), x+y=y+x,
\\ x+0=x, x+(-x)=0,-(x+y)=-x-y,
\end{array}} \right\}
\end{split}
\end{equation}
\begin{equation}
\label{4.5.2} \gamma^{-1}\gamma x = x, 0\gamma= 0, (x + y)\gamma =
x\gamma + y \gamma, \gamma \in \Delta_0.
\end{equation}

Besides that we  will  use the free group $G=<\alpha, \beta
>$ with the basis ${\{}\alpha, \beta {\}}$.
Let's remember that every element of  group $G$ is presented by
the uniquely  introduced word, that is a word which does not have
sub-words  of type  $\gamma \gamma^{-1}, \gamma\in\Delta_0$.  At
the same time it  is supposed that $(\gamma^{-1})^{-1}=\gamma$.

{\bf Lemma.} {\it In variety algebras $U(\Delta_1, S)$ the problem
of the equality of words for the free algebras is solved.}

\textbf{Proof of Lemma.} Let $F(A)$  be a   free algebra of
variety $U(\Delta _1 ,S)$ with the basis $A$. Lets introduce the
meaning of the canonical $\Delta_1$-word in alphabet $A$: the
$\Delta_1$-word $R$ in alphabet $A$ we call  canonical if $R=0$ or
has the form
\begin{equation}
\label{4.5.3}
R=c_{11}a_{1}\gamma_{11}+\ldots+ c_{1k_{1}} a_{1}\gamma_{1k_1} +
\ldots + c_{nk_1}a_{n}\gamma_{nk_1} +\ldots + c_{nk_{n}}
a_{n}\gamma_{nk_n},
\end{equation}
where $n>0$, $c_{ij} \in Z \backslash {\{}0{\}}$, for $k_{i} \ge
0$, $\gamma_{ij} \in G$; and $\gamma_{ij}$, for $j=1,\ldots,
k_{i}$, are various  pairs of  introduced words  from the group
$G$ for any fixed $i \in {\{}1,\ldots,n{\}}$.

Here under $сa_{i}\gamma$ it is necessary  to understand  the sum
of  $c$ items $a_{i}\gamma$ for  $c>0$ and $(-c)$ items for $c<0$.
The items of type $c_{1} a_{i}\gamma$, $c_{2}a_{i}\gamma$ will be
called similar  and the changing of its sum by the words
$(c_{1}+c_{2})a_{i}\gamma$ will be called \emph{introducing of
similar}. It is evident that the \emph{introducing of similar}
should be realized by the elementary transformations of words on
identities (\ref{4.5.1}), (\ref{4.5.2}).

Let's note that in sum  (\ref{4.5.3})  brackets which determine
the order of operations were not used. Anyway the equality of
(\ref{4.5.3}) is correct due to the   law of associativity of
addition. The law of associativity will be used with this meaning
in future.

Let's prove that for any $\Delta_1$-word $P$ in alphabet $A$ there
exists the equivalent of its canonical $\Delta_1$-word $R$ and
this word is unique within all permutations of  items.

The existence of the word $R$ we prove  by  induction with  the
rank of the word $P$.

If  $rankP=0$, than $P=0$ or $P=a_{i}$  and the assertion is
obvious. Let's assume that it is true for all words of rank $r<m$
and we will consider the case where  $rankP=m>0$.

According to the definition of the $\Delta_1$-word it is possible
is three cases:

1) $P=(P_{1})+(P_{2})$. In this case to find the  unknown word $R$
it is enough to find the sum of the canonical words for $P_{1}$,
$P_{2}$, and  in it, using  identities (\ref{4.5.1}), to group
items by the same factors from $A$ in each resulting sum and
reduce similar items  and delete  zero items, if they exist.

2) $P=(P_{1})\gamma, \gamma \in \Delta_0.$ In this case it is
enough using identities (\ref{4.5.2}) to add  to every item of the
canonical word for $P$ the letter $\gamma$ and make a reduction of
type $\gamma^{-1}\gamma$, if it exists.

3) $P=-(P_{1})$. In this case it is  enough in the canonical word
for $P_{1}$ to change all the coefficients from $Z$ to  the
opposite numbers.

From this result the existence of the canonical word equivalent to
the  word $P$ is proved.

Now we will prove its uniqueness, within the permutations of
items. Firstly,  from the demonstrated proof of existence, we will
obtain the procedure of introducing   any given $\Delta_1$-word
$P$ to a canonical type.

Step 1. Using identities
$$
(x+y)+z=x+(y+z), -(x+y)=-x-y,\,\,\, (x + y)\gamma = x\gamma +
y\gamma, \gamma \in \Delta_0,
$$

we open all brackets in word $P$ and bring it to the sum of words
of type $c_{i}a_{i}\gamma_i$.

Step 2. Using identities  $\gamma^{-1}\gamma x = x, \gamma \in
\Delta_0$,  in every word of type $c_{i}a_{i}\gamma_i$ change the
word $\gamma_i$ from the group $G$ by using the  introduced word.

Step 3. In the received sum  reduce  similar words.

Step 4. Using the identities $0+x=x$, $x+0=x$, we delete the zero
items.

Step 5. We group the items with factors from $A$ in the remaining
items.

We notice that the described procedure  is  not an algorithm due
to  its ambiguity  which may appear in steps 3, 5. Therefore the
resulting canonical words  should be different but it will be easy
to see, that they  differ only in the re-arrangement of their
items.

The set of all such received  canonical words will be  denoted by
$\Phi(P)$. In this way, all words  from $\Phi(P)$ are canonical
and equivalent  to the word $R$ and they differ  from each  other
only by the re-arrangement of  items in sums with the same factors
from $A$. In  this case such a permutation may be arbitrary.

Now we will prove that any canonical word which is equivalent to
$P$ is contained in $\Phi(P)$. Let $R_{1} \in \Phi(P)$ and $R_{2}$
be  any canonical word equivalent to $P$ (which is received,
possibly  by  another method). Then the words $R_{1}$, $R_{2}$ are
equivalent and so  it is possible to pass from $R_{1}$ to $R_{2}$
with the help of a finite number of elementary transformations:
$$
R_{1 } \to T_{1} \to T_{2} \to  \ldots  \to T_{m}= R_{2}.
$$

We next prove that $\Phi(T_{i})=\Phi(T_{i + 1})$, $i=1, \ldots, m
-1.$ For this it is  necessary to consider all possible elementary
transformations $g:T_{i} \to T_{i + 1}$. At the same time this
transformation by the  identity of  associativity should not be
considered, because  it has been used implicitly in the recording
of sums of many items without the re-arrangement of brackets.

a) Let the transformation $g:T_{i} \to T_{i + 1}$  be carried out
by the  identity $x+y=y+x$. Then the sums received from words
$T_{i}$, $T_{i + 1},$ following  step 1   of procedure $\Phi$ will
differ only by  the order of items so  after steps 2-3 we receive
words, which are differed only in the order of items. Hence in
this case $\Phi(T_{i})=\Phi(T_{i + 1})$.

b) Let  the transformation $g:T_{i} \to T_{i + 1}$  be  carried
out by the identity $-(x+y)=-x-y$. Then using  the words $T_{i}$,
$T_{i + 1}$ by the procedure  $\Phi,$ after the first step we get
the same words   and  hence $\Phi(T_{i})=\Phi(T_{i + 1})$.

c) Let  the transformation $g:T_{i} \to T_{i + 1}$  be carried out
by the identity $\gamma^{-1}\gamma x = x$. Then, using the  words
$T_{i}$, $T_{i + 1}$ by the procedure $\Phi,$ after the first step
we  accordingly get the sums $S_1, S_2 $, moreover, some items in
$S_2$ will be obtained from the corresponding  items of  the sum
$S_{1}$ by deleting  $\gamma^{-1}\gamma$, and so after  step 2 we
will obtain the same sums. Hence  in this case we will have the
equality $\Phi(T_{i})=\Phi(T_{i + 1})$.

This equality  by analogy proves those cases where $g$ carried out
on other identities from (\ref{4.5.1}), (\ref{4.5.2}).

From the proved equalities  $\Phi(T_{i})=\Phi(T_{i + 1})$,
$i=1,\ldots ,m-1$ it follows that
$$
\Phi(R_{1})=\Phi(R_{2}).
$$

This way all the canonical words which are equivalent to word $P$
differ only by the  re-arrangement of items in sums with the same
factors from $A$. From this  we have received the proof of the
lemma.

For recognition of the equivalence of the words $P_{1}$, $P_{2}$
it is enough to find and compare the sets $\Phi(P_{1})$,
$\Phi(P_{2})$. $P_{1}$ and $P_{2}$ are equivalents in the algebra
$F(A)$ if and only if $\Phi(P_{1})= \Phi(P_{2})$. Evidently  for
testing of the equality $\Phi(P_{1})=\Phi(P_{2})$ it is enough to
compare only  one representative from $\Phi(P_{1})$ and
$\Phi(P_{2})$.

Let's return now to the proof of   theorem 1.

Following [13,14,18]  let's establish the connection between the
varieties $Q(\Omega, \Sigma)$ of all T-quasigroups and the variety
of algebras $U(\Delta_1, S)$. Since the signature  $\Delta_1$
contains the 0-ary operation $O$, then  to establish the
equivalence of the varieties it is necessary to expand $\Omega$ by
way of introducing the 0-ary operation $O$. In this connection, we
will be considering the variety $Q(\Omega_1, \Sigma)$, where
$\Omega_{1}={\{}\cdot, /, \backslash , u{\}}$,  and $u$   is a
symbol of 0-ary operation.

It is known (and evident) that if a quasigroup is isotopic to
group $G$ then it is principally isotopic to some other group
which is isomorphic to $G$. Since the variety of groups is closed
compared   to isomorphisms of the groups, then   it should be
clear that the considered quasigroups are principally isotopic to
groups from $U$. So in contrast with [13,14], we will not
introduce in signature $\Delta_1$ the additional symbol of unary
operation $c$. The transfers of signature $\tau :\Omega _1 \to
\Delta _1 $ and $\sigma :\Delta_1 \to \Omega_1$ will be determined
by on the following formulas:
$$
\tau (\cdot )(x,y)=x\alpha + y\beta,
$$
$$
\tau (/)(x, y)=(x - y\beta)\alpha^{-1},
$$
$$
\tau (\backslash)(x, y) = (-x\alpha + y)\beta^{-1},
$$
$$
\tau(u)=0,
$$
$$
\sigma (+)(x, y)=(x/u)(u/u)\backslash y),
$$
$$
\sigma (-)(x)=(u/u)((x/u)\backslash u),
$$
$$
\sigma (\alpha )=x(u\backslash u), \sigma (\beta )=(u/u)x,
$$
$$
\sigma(\alpha^{-1}) = x/(u\backslash u), \sigma
(\beta^{-1})=(u/u)\backslash x,
$$
$$
\sigma (0)=u.
$$

The transfers  $\tau $  and  $\sigma$  induce  one-to-one maps
between $\Omega_1$-words and  $\Delta_1$-words and vice versa. For
the $\Omega_1$-word $P$, the  $\Delta_1$-word $\tau(P)$  is
obtained by the use of all operations from $\Omega_1$ via its
$\tau$-transfers. Analogously  for the $\Delta_1$-word $R$ the
$\Omega_1$-word $\sigma(R)$ is determined. Let's note that
according to the result of  F.N.Sokhatsky [8], the system of
identities of varieties of quasigroups $Q(\Omega_1, \Sigma)$ which
are isotopic to  abelian groups is obtained by addition to the
system of identities $\Sigma_0$ which are received from
(\ref{4.5.1}) by changing the  group operation  by its transfers.

From [13,14] it follows that the transfers $\tau$ and $\sigma$
also  induce  the bijective and one-to-one maps $T_\tau, T_\sigma$
of the algebras of varieties $U(\Delta_1, S)$ and $Q(\Omega_1,
\Sigma)$ and the above mentioned varieties are rationally
equivalent. Hence, two $\Delta_1$-words are equivalent in the free
algebra $F(А)$ of variety $U(\Delta_1, S)$ if and only if its
transfers are equivalent in free quasigroup $Q(A_{1})$ with the
basis $A_{1}=A \cup {\{}u{\}}$ of the  variety $Q(\Omega_1,
\Sigma)$. Evidently,  the classes of words which form elements of
the free quasigroup $Q(A)$ of variety $Q(\Omega, \Sigma)$ are
subclasses of some classes of quasigroup $Q(A_{1})$. Hence the
elements from $Q(A)$ are equal if and only if  the corresponding
elements in algebra $F(A)$ are equal. From here and the above
lemma follows the proof of theorem 1.

{\bf Corollary 1.1}. {\it In a variety of  all medial quasigroups
the word problem  for  free algebras is solved}.

{\bf Remark}. There is also a more general approach to the concept
of a linear quasigroup, namely,  quasigroups which are linear over
some loop are considered (T. Kepka [22,23], P. Nemec [24], V.A.
Shcherbacov [25] and others). A quasigroup $(Q, \cdot)$  is {\it
called linear over loop} $(Q, +)$, if it has the form $xy =
\left(\varphi x + \psi y\right) + d$, where $\varphi, \psi \in
Aut\left(Q, + \right), d \in Q,$ assuming that  the loop $(Q, +)$
will be fairly well-known and examined loop, for example Moufang
loops, i.e. loops with the identity $x+(y+(x+z))=((x+y)+x)+z$ for
all $x,y,z \in Q$.

The general idea of a quasigroup is linear over  loop has been
crystallized  in the works of Prague's algebraic school (T.Kepka,
J.Jezek, [18,22-24]). Recently in  literature   the term {\it
generalized linear quasigroups} [25] has appeared.

As V.A. Shcherbacov [25] noted, many well known (classical)
objects belonged to  the class of generalized linear quasigroups.
For example, medial quasigroups ({\it $\text{Toyoda Theorem} \,
[26]$}), distributive quasigroups ({\it Belousov Theorem} [27],
distributive Steiner quasigroups, leftdistributive quasigroups
({\it Belousov-Onoi Theorem} [28], CH-quasigroups ({\it Manin
Theorem} [29], T-quasigroups ({\it Belyavskaya Theorem} [30],
$n$-ary groups {\it(Gluskin-Hosszu Theorem} [31], $n$-ary medial
quasigroups ({\it Evans Theorem} [32] {\it and Belousov Theorem}),
F-quasigroups ({\it Kepka-Kinyon-Phillips Theorem}, [33]  are
quasigroups of such kind. It is necessary to note also  the new
results of V.A. Shcherbacov on the structure of left and  right
F-, SM- and E-quasigroups ({\it Shcherbacov Theorem)} [34].

The generalization was suggested in view  of  several theorems
about connections between some classes of quasigroups and loops.
The first in this series is Toyoda-Mudoch Theorem  on medial
quasigroups. Any medial quasigroup can be obtained as follows: $xy
= \varphi x + \psi y + d$, where  $\varphi, \psi \in Aut(Q, + )$,
$\varphi \psi = \psi \varphi $, $d \in Q$, $(Q, + )$ is abelian
group. A quasigroup $(Q, \cdot)$ with identities   $x \cdot yz =
xy \cdot xz, \quad xy \cdot z = xz \cdot yz$ is called \emph
{distributive}. If a quasigroup satisfies only the first identity,
then it is called \emph {leftdistributive}. In 1958 V.D. Belousov
proved that every distributive quasigroup can be obtained as
follows: $xy = \varphi x + \psi y$, where $\varphi$ and $\psi$
some automorphisms commutative Moufang loop (CML) $(Q, + )$.

A quasigroup with identities  $ xy = yx, \quad x (xy) = y $ is
called \emph{CH-quasigroup} such that in it  any three elements
generate medial subquasigroup. CH-quasigroups were introduced by
Yu.I. Manin in connection with solving the problem in algebraic
geometry, namely researching  cubic hypersurfaces. Yu.I. Manin
proved that any CH-quasigroup can be obtained by using the
following construction:  $xy = (-x - y) + d$, where $d \in Z(Q, +
)$ and $Z(Q, + ) = \{a \in Q\vert a + (x + y) = (a + x) + y,$ $x,
y \in Q\}$ is the center of CML. Later study of linear quasigroups
over Moufang loops, CML, groups, abelian groups was also carried
out by other mathematicians.

By a known  Stein Theorem [35] any leftdistributive quasigroup
$(Q,\cdot)$ which is  isotopic to the group $(Q,+)$, can be
obtained using the following construction: $xy = x + \varphi (-x +
y),$ where  $\varphi$ is automorphism of $(Q,+)$. Due to the
associativity of the group operation, we obtain: $xy = (x-\varphi
x) + \varphi y = \psi x + \varphi y$. An automorphism $\varphi$
such that  $\psi $ is substitution. Thus, a leftdistributive
quasigroups, isotopic to groups, in fact, is  right linear over
groups.

After well know  work of V.D.Belousov  [36] czech algebraists --
T. Kepka, P. Nemec, J. Jezek and representatives Belousov's
quasigroup school  - G.B. Belyavskaya, V.A. Shcherbacov, V.I.
Izbash, K.K. Shchukin, F.N. Sokhatsky, P.N. Sirbu, A. KH. Tabarov,
W.A. Dudek  comprehensively and intensively  studied a linear
quasigroups and some of their generalizations. They investigated
the algebraic (morphisms, congruences, core, center, associator,
commutator, multiplication group) and the combinatorial
(orthogonality, numerical estimates, Latin squares) aspects of
generalized linear quasigroups , also $n$-ary linear quasigroups.

It is important to note that  Belousov's quasigroup school has
become a world center for the development of the theory of
quasigroups and loops. In addition, in the development of the
theory of quasigroups and loops enormous contribution have been
made by the representatives of the Czech  algebraic school.
Sufficiently detailed historical overview of the theory of
quasigroups is contained in doctoral dissertations of H. Kiechle
[37] and V.A. Shcherbacov [25].

\begin{center}
{\bf II. About one approach  of finding quasigroup identities from
some variety of loops}
\end{center}

\quad This part  is devoted to an  approach  of finding the
identities in the class of quasigroups isotopic to known classes
of loops from some  variety of loops. It is assumed  that the
class of loops from  some  variety of loops are given by an
identity or system of identities. For  this,  the notion  a
derived identity is introduced and it is proved that for any
identity  from  some of a variety  of quasigroups (loops)  there
is derivative identity. The introduced notion of  derived
identities   allows to find an arbitrary identity for the class of
quasigroups  which are isotopic not only to  groups but also to
loops from some variety of  loops and generalizes the method of
A.A. Gvaramiya [38], where for the class of quasigroups are
isotopic to a groups, it is possible  to obtain  any quasigroup
identity from the  group identities. Also by this way it is easy
to obtain V.D. Belousov's identities which  characterized  the
class of quasigroups is isotopic to a group (an abelian
group)[36]. As an illustration we show that the class of
quasigroups which is isotopic to the nilpotent group is
characterized by identity. It should be noted that the notion of a
derived identity in terms of free quasigroup and the theory of
automats  also can be found  in [38]. However, our proposed
approach does not require using  free objects, in particular a
free quasigroups. Sufficiently confined by the methods of  theory
of quasigroups. It should be noted that henceforth we shall
consider a quasigroup as an algebra $(Q, \cdot, /,  \setminus )$
with three binary operations $( \cdot )$, $( / )$ and $( \setminus
)$, satisfying the following identities:
$$
(xy)/y=x, \,\,(x/y)y=x,\,\, y(y \backslash x)=x,\,\, y \backslash
(yx)=x.
$$

{\bf Definition 2.1.} {\it Let $\left(Q,\cdot \right)$ be a
quasigroup  with the identity $w_1 = w_2 $. The identity $w'_1 =
w'_2 $ is called  derived from  the identity $w_1 = w_2 $, if it
is received from $w_1 = w_2$ by adding two  new variables $u$, $v$
in the left and the right sides of this identity using the
quasigroup operations $( \cdot ), (\backslash )$, $(/)$.}

In this case the identity $w_1 = w_2 $  is considered to be  main.
We  give  some examples of   derived  identities: the identity $(x
\backslash z ) \backslash( y  / t )= (x / z )/ (y / t)$ is
received from the identity  $x \backslash z = x / z$ by adding two
new variables $y $, $t$ and the operations ($\backslash )$, (/).
In some  cases the operations $(\cdot )$, ($\backslash )$, (/) may
be coincided. In this case we have the identity from one
operation, for example the identity $(xt \cdot y) \cdot uv = xu
\cdot (y \cdot zv)$ which is received from the identity $(xy)
\cdot (uv) = (xu) \cdot (yv)$ (the identity of mediality) by
replacing the variables $x$ on $xt$, $v$ on $zv$.

As in the case of balanced identities (see [36]), the derived
identity $w'_1 = w'_2$ is called {\it the first kind} if variables
are ordered identically and {\it the second kind}  otherwise. For
example, the  identity $(x \backslash z ) \backslash( y  / t )= (x
/ z )/ (y / t)$  is the derived identity of first kind, the
identity $(xy) \cdot (uv) = (xu) \cdot (yv)$ is the derived
identity of second kind.

{\bf Lemma 2.1.}  {\it For any  identity $w_1 = w_2 $ of some
variety quasigroups there is a derived identity $w'_1 = w'_2$.}

{\bf Proof:} Let $\left(Q,\cdot \right)$ be a quasigroup. By the
Albert's theorem [1] the quasigroup $\left(Q,\cdot \right)$ is
isotopic to some loop $\left( {Q,\mbox{ } + } \right)$. It is
sufficiently to consider the principal isotopy $x+y = R_{a}^{-1} x
\cdot L_{b}^{-1} y$. Suppose that in the loop is satisfied the
identity $w_1 = w_2 $. Now we transfer  the quasigroup using the
isotopy $T=(R_a, L_b, \varepsilon )$. We obtain an identity
containing the permutations $ R_a $ and $ L_b $ (or $R_{a}^{-1}$
and $L_{b}^{-1})$. Since the elements $a,b \in Q$ may run over all
$Q$, so replacing $a$ and $b$ on any $u,v \in Q,$ we have the
identity $w'_1 = w'_2$ which is obtained from the identity $w_1 =
w_2 $ by adding two new variables $u$ and $v$.

\endproof

{\bf Remark 2.1}. It is obvious that forming an isotopy we can
obtain various derived identities. So for  some identity there are
indefinite   derived identities, in the case when  a quasigroup is
an indefinite.

{\bf Remark 2.2}. It is possible to obtain also such is called an
identity with permutations for some identity [39].

{\bf Example.} Let $(Q,\cdot)$  be a quasigroup  isotopic to the
Moufang loop $(Q, + )$:  $x + y = R_{a}^{ - 1} x \cdot L_{b}^{-1}
y$. A loop $(Q, + )$ is called Moufang loop if in it is satisfied
the identity:
$$
(x + (y + z)) + x = (x + y) + (z + x),
$$
for any $x,y,z \in Q$.  At  transition to quasigroup operation and
replacing the elements $a,b \in Q$ by  $u, v \in Q$ we have the
following identity:
$$
((x \cdot (v \backslash ((y \cdot (u\backslash zu))) / u) \cdot
(u\backslash xu) = (((x \cdot (u\backslash yu)) / u)) \cdot
(v\backslash (z \cdot (u\backslash xu))).
$$

{\bf Lemma 2.2.}  \emph{For  any  derived identity  $w'_1 = w'_2 $
there is a main identity $w_1 = w_2 $ and vice versa.}

Denote by $\left| w \right|$ the length of an identity, that is
number of variables in it.

{\bf Lemma 2.3.} Suppose that for a quasigroup $(Q,\cdot)$ an
identity  $w_1 = w_2 $ is the main and the identity $w'_1 = w'_2$
is its  derived identity,  $\left| {w_1 } \right| = m$  and
$\left| {w_2 } \right| = n$. Then $\left| w \right|=m + n$ and the
length of derived identity is  equal  to $k=\left| {w'}
\right|=\left| {w'_1 } \right|+\left| {w'_2 } \right| \le m +
n+2$.

The proofs of lemmas 2.2 and 2.3 are obvious.

Now let us show how  to obtain the identities which characterize
the class of quasigroups  isotopic to the nilpotent group of level
$ \le $n. This class  is described by identity    (4) (see Theorem
2.1).
The corollary of this result is theorem V.D. Belousov [36], where
quasigroups which are  isotopic to an abelian groups are
characterized by the identity. In [36] this identity are given in
other form and by the operations $\left( \cdot \right)$  and
$\left( \backslash \right)$, namely
\begin{equation}
\label{1} x\backslash (y(u\backslash v)) = u\backslash
(y(x\backslash v)).
\end{equation}

Let us show that in case  when a quasigroup is isotopic to an
abelian group we can  obtain other  identity  from Belousov's
identity using transformations. Indeed, we can re-write the
identity (\ref{1}) in the following form:
$$
L_{x}^{-1} (y \cdot L_{u}^{-1} v) = L_{u}^{-1} (y \cdot L_{x}^{-1}
v)
$$
$$
L_{u} L_{x}^{-1} (y \cdot L_{u}^{-1} v) = (y \cdot L_{x}^{-1} v).
$$

Now $v$ can be replaced by  $L_{x}^{-1} v$. Then
$$
L_{u} L_{x}^{-1} (y \cdot L_{u}^{-1} L_{x}v) = (y \cdot v),
$$
$$
y \cdot L_{u}^{-1} L_{x} v = L_{x} L_{u}^{-1} (y \cdot v),
$$
$$
y \cdot L_{u}^{-1} (x \cdot v) = x \cdot L_{u}^{-1} (y \cdot v),
$$
$$
y \cdot L_{u}^{-1} (R_{v} x) = x \cdot L_{u}^{-1} (R_{v}y).
$$

In the last equality  $x$ and $y$ can be replaced by $R_{v}^{-1}
x$ and $R_{v}^{-1} y$ accordingly. Then  $R_{v}^{-1} y \cdot
L_{u}^{-1} x = R_{v}^{-1} x \cdot L_{u}^{-1} y, \quad (y / v)
\cdot (u\backslash x) = (x / v) \cdot (u\backslash y)$, or
$$
(x / v) \cdot (u\backslash y) = (y / v) \cdot (u \backslash x).
$$

Below we show that the last identity is (\ref{3}) when $n=1$.

Also note that  quasigroups are isotopic to  abelian groups can be
characterized by other identity  which is "symmetrical" to
Belousov's identity [36] which includes the operations $(\cdot)$
and $(/)$, namely
$$
((u / v) x) / y = ((u / y) x) / v.
$$

It is easy to show that from the last identity we can  obtaine
Belousov's identity by such transformations.

Now for convenience we will use the additive form  of a group. As
well known [40] a commutator of the elements $x_{1}$ and $x_{2}$
in a group $(Q, +)$ is $\left[{x_{1}, x_{2}}\right] = x_1 + x_2 -
x_1 - x_2 $. A commutator of the elements $x_1, x_2,...,x_{n + 1}
$ is defined  recursively:
$$
\left[x_1, x_2,..., x_{n + 1}\right] = \left[{...}\right.\left[
{x_1, x_2} \right], x_3 \left. \right],...,x_{n + 1}\left.
\right].
$$

A group $\left( {Q,\mbox{ } + } \right)$ is called  nilpotent of
class  $ \le $n, if
$$
\left[ x_1, x_2,...,x_{n + 1}\right] = 0,
$$
where 0 is zero element of the group $(Q, +)$.

Let us denote $\left[x_1, x_2 \right] = x_1 + x_2 - x_1 - x_2 =
t$, where $t \in Q$. Then
\begin{equation}
\label{2} x_1 + x_2 = t + (x_2 + x_1).
\end{equation}

From $xy = R_{a} x + L_{b} y$  it follows that $x + y = R_{a}^{-1}
x \cdot L_{b}^{-1} y$. Hence we can  write (\ref{2}) in the
following form:
$$
R_{a}^{-1} x_{1} \cdot L_{b}^{-1} x_{2} = R_{a}^{-1} t \cdot
L_{b}^{-1} (R_{a}^{-1} x_2 \cdot L_{b}^{-1} x_{1}).
$$
Hence
$$
R_{a}^{-1} t = (R_{a}^{-1} x_{1} \cdot L_{b}^{-1} x_{2}) /
L_{b}^{-1} (R_{a}^{-1} x_{2} \cdot L_{b}^{-1} x_{1})
$$
or
$$
t = R_{a} ((R_{a}^{-1} x_1 \cdot L_b^{-1} x_2) / L_{b}^{-1}
(R_{a}^{-1} x_2 \cdot L_{b}^{-1} x_{1})).
$$

So, to  commutator $\left[x_1, x_2\right]= x_{1} + x_{2} - x_{1} -
x_{2} $ of the elements $x_1 $ and $x_2 $ of a group $(Q, +)$
corresponds  the  element $t = R_a ((R_a^{-1} x_1 \cdot L_{b}^{-1}
x_2 ) / L_{b}^{-1} (R_{a}^{-1} x_2 \cdot L_b^{-1} x_1 ))$ of a
quasigroup $(Q,\cdot)$.

For convenience   denote    the element $t$ by $\left\{x_{1},
x_{2}\right\}$ and  is called  a quasicommutator of the elements
$x_{1} $ and  $x_{2}$ in a quasigroup $\left(Q, \cdot \right)$:
$$
\left\{ x_1, x_2\right\} = R_{a} ((R_a^{-1} x_{1} \cdot L_b^{-1}
x_2 ) / L_b^{-1} (R_a^{-1} x_{2} \cdot L_{b}^{-1} x_1)).
$$

Similarly to  commutator $\left[\left[ x_1, x_2 \right], x_3
\right]$ of the elements $x_1 $, $x_2 $ and  $x_3$ of a group $(Q,
+)$ corresponds  the  element $R_a ((R_a^{-1} \left\{x_1,
x_2\right\} \cdot L_b^{-1} x_3 ) / L_b^{-1} (R_a^{-1} x_3 \cdot
L_b^{-1} \left\{ x_1, x_2 \right\}))$ which is denote by
$\left\{x_1, x_2, x_3\right\}$ and is   called a quasicommutator
of the elements $x_1 $, $x_2 $, $x_3$. So a quasicommutator of the
elements $x_1, x_2,...,x_{n + 1}$ is defined recursively:
$$
\{x_1, x_2,...,x_{n + 1}\}=
$$
$$
=R_a ((R_a^{-1}\left\{ x_1, x_2,...,x_n \right\} \cdot L_b^{-1}
x_{n + 1}) /L_b^{-1} (R_{a}^{-1} x_{n + 1} \cdot L_{b}^{-1}
\left\{x_{1}, x_{2},...,x_n \right\})).
$$

{\bf Theorem 2.1.} {\it A quasigroup $(Q,\cdot)$ is isotopic to a
nilpotent group of class $\le n$ if and only if in it the
following identity is satisfied:}
\begin{equation}
\label{3} \left(\left\{x_1, x_2,..., x_n \right\} / u\right) \cdot
(v \backslash x_{n + 1}) = (x_{n + 1}/u) \cdot (v \backslash
\left\{ x_1, x_2,..., x_n \right\}).
\end{equation}

{\bf Proof:} Let $(Q,\cdot)$ be  quasigroup is isotopic to a
nilpotent group of class $\le n$. We show that in  quasigroup
$(Q,\cdot)$ the identity (\ref{3}) is satisfied. As was noted
above it is enough to consider a principle isotopy of the form $x
+ y = R_{a}^{-1} x \cdot L_{b}^{-1} y$. The proof can be shown by
induction.

If $n=1$, then  the group is  abelian and identity (\ref{3}) has
the form
\begin{equation}
\label{4} (x_1 / u) \cdot (v\backslash x_2 ) = (x_2 / u) \cdot
(v\backslash x_{1}).
\end{equation}

Let us show that in this case (\ref{4}) is satisfied. For this
first note that from $xy = R_{a} x + L_{b} y$ follows $x / y =
R_{a}^{-1} (x + JL_{b} y)$ and  $x \backslash y = L_{b}^{-1}
(JR_{a} x + y)$, where  $Jx = - x$. Then we have:
$$
(x_1 / u) \cdot (v \backslash x_2 ) = R_a R_a^{ - 1} (x_1 + JL_b
u) + L_b L_b^{ - 1} (JR_a v + x_2 ) = x_1 + JL_b u + JR_a v + x_2,
$$
$$
(x_2 / u) \cdot (v\backslash x_1 ) = R_a R_a^{ - 1} (x_2 + JL_b u)
+ L_b L_b^{ - 1} (JR_a v + x_1 ) = x_2 + JL_b u + JR_a v + x_1 .
$$

Comparing the right sides of last  equalities and  accordingly
commutativity of a group  concluded  that the equality (\ref{4})
is satisfied.

Now by induction suppose that the identity (\ref{3}) in the
quasigroup $(Q,\cdot)$ is satisfied which is isotopic to a
nilpotent group $(Q, +)$ of class  $n=k$, i.e.
$$
(\left\{x_1, x_2,...,x_{k - 1}\right\} / u) \cdot (v\backslash
x_k) = (x_k / u) \cdot (v\backslash \left\{ x_1,
x_2,...,x_{k-1}\right\}).
$$

Let us show the case when $n=k+1$. In this case the identity
(\ref{3}) has the form:
$$
(\left\{ x_1, x_2,..., x_k \right\} / u) \cdot (v \backslash x_{k
+ 1}) = (x_{k + 1} / u) \cdot (v \backslash \left\{x_1,
x_2,...,x_k \right\}).
$$

Replacing the operation $(\cdot)$ by $(+)$ in the left and right
sides of the last  identity we obtain:
$$
(\left\{ x_1, x_2,..., x_{k}\right\} / u) \cdot (v \backslash x_{k
+ 1} ) = R_{a} R_{a}^{-1} (\left\{ x_1, x_2,...,x_k\right\}+
$$
$$
+JL_{b} u) + L_b L_b^{-1} (JR_a v + x_{k + 1})=\left\{x_1,
x_2,...,x_k\right\} + JL_b u + JR_{a} v + x_{k + 1},
$$
$$
(x_{k + 1} / u) \cdot (v \backslash \left\{ x_1, x_2,...,x_k
\right\}) = R_a R_a^{-1} (x_{k + 1} + JL_b u) + L_b L_b^{-1} (JR_a
v +
$$
$$
+\left\{x_1, x_2,...,x_k\right\})= x_{k + 1} + JL_b u + JR_{a} v +
\left\{ x_1, x_2,...,x_k\right\}.
$$

On other side in a nilpotent group of class  $n=k$:
$$
\left[x_1, x_2,...,x_{k + 1}\right] = [\left[x_1, x_2,...,х_k
\right], x_{k + 1} ] = 0
$$
or
$$
[x_1, x_2,...,x_k] + x_{k + 1} = x_{k + 1} + [x_1, x_2,...,x_k ].
$$

But the commutator $[x_1, x_2,...,x_k]$ corresponds to the
quasicommutator $\left\{x_1, x_2,...,x_k\right\}$. Then
$$
\left\{x_1, x_2,..., x_k \right\} + x_{k + 1} = x_{k + 1}
 + \left\{x_1, x_2,...,x_k\right\}.
$$

As the element $x_{k +1}$ is  arbitrary, then suppose $x_{k + 1} =
JL_b u + JR_a v$ we get that the element $JL_b u + JR_{a} v$
commutes with the quasicommutator $\left\{x_1,
x_2,...,x_k\right\}$:
$$
\left\{ x_1, x_2,...,x_k\right\} + JL_b u + JR_{a} v = JL_b u +
JR_a v + \left\{x_{1}, x_2,...,x_k\right\}.
$$
Then we have
$$
\left\{ x_1, x_2 ,...,x_k \right\} + JL_b u + JR_a v + x_{k + 1} =
x_{k + 1} + JL_b u + JR_{a} v + \left\{x_1, x_2,...,x_k\right\}.
$$
The last equality shows that the identity (\ref{3}) is true for
$n=k+1.$

Now let in a quasigroup $(Q,\cdot)$  the identity (\ref{3}) is
satisfied
$$
(\left\{x_1, x_2,...,x_{n}\right\} / u) \cdot (v \backslash x_{n +
1}) = (х_{n + 1} / u) \cdot (v\backslash \left\{{x_1, x_2,...,x_n}
\right\}).
$$

It is necessary to prove that a quasigroup is isotopic to a
nilpotent group of class  $\le n$.

By induction for  $n=1$ have:
$$
(x_1 / u) \cdot (v \backslash x_2 ) = (x_2 / u) \cdot (v
\backslash x_1 )
$$

As we have shown in above the last identity  characterizes the
class of quasigroups which are  isotopic to  abelian groups.

For $n=2$ we have the following identity:
$$
(\left\{x_1, x_2\right\} / u) \cdot (v \backslash x_3 ) = (x_3 /
u) \cdot (v \backslash \left\{x_1, x_2\right\}),
$$
which  characterizes a quasigroups  isotopic to a metabelian
groups (a nilpotent group of class  2). Continuing this process,
we find that a quasigroup is isotopic to a nilpotent group of
class $ \le n.$

\endproof

{\bf Corollary 2.1.} {A quasigroup $\left( Q,\cdot\right)$ is
isotopic to a nilpotent group of class $\le n$ if and only if in
it  the following equality is satisfied:
$$
\{x_1, x_2,...,x_{n + 1}\}=ba,
$$
where $ba = 0$ is zero element of the group $\left(Q,+ \right)$}.

A quasigroup is isotopic to a commutative group, can also be
characterized  by  an identity with  4 variables, namely

{\bf Corollary 2.2. (Belousov's Theorem [2])}  {\it A quasigroup
$\left( Q,\cdot\right)$ is isotopic to an abelian group  if and
only if in it the following equality is satisfied:}
$$
 x\backslash (y(u\backslash v)) = u\backslash (y(x\backslash v)).
$$

{\bf Corollary 2.3.} {\it A quasigroup $\left( Q,\cdot\right)$ is
isotopic to an abelian group  if and only if in it  the following
equality is satisfied:}
\begin{equation}
\label{5} \left( {\left( {u / v} \right)  x} \right) / y = \left(
{\left( {u / y} \right)  x} \right) / v.
\end{equation}

The variety of  quasigroups isotopic to  abelian groups   was also
investigated by M.M. Glukhov,  A. Drapal,   etc. In particular,
M.M. Glukhov has described  an identities of length 4, which
characterizes class of   quasigroups which are isotopic to abelian
groups:
\begin{equation}
\label{6}((x/y)/u)/v = ((x/v)/u)/y,
\end{equation}
\begin{equation}
\label{7}(x(y\backslash(uv)) = u(y\backslash(xv))
\end{equation}
A. Drapal [41] shows that this class can be characterized by the
following identity:
\begin{equation}
\label{8}((xy)/u)v = ((xv)/u)y
\end{equation}

Note also that independently from others  M.M. Glukhov has
obtained the identities  (7), (10-12).

As known [40] a group is called {\it engel}  if in it the
following identity is satisfied:
$$
[x, \underbrace {y,...,y}_{n\ \text{copies}}] = 0.
$$

Using Theorem 2.1 we can  get criterion when a quasigroup is
isotopic to  engel group.

{\bf Corollary 2.4.} {\it A quasigroup $(Q,\cdot)$ is isotopic to
an engel group if and only if in it  the identity:}
\begin{equation}
\label{9} (\{x, \underbrace{y,...,y}_{(n - 1)\ \text{copies}}\} /
u) \cdot (v \backslash y) = (y / u) \cdot (v \backslash \{ x,
\underbrace{y,...,y}_{(n - 1)\ \text{copies}} \}).
\end{equation}
is satisfied.

{\bf Remark 2.3}. Quasigroups are isotopic to  abelian groups also
can be characterized by the identities which  contain 5 variables.

{\bf Theorem 2.2.} {\it A quasigroup $\left( Q,\cdot\right)$ is
isotopic to an abelian group  if and only if in it  at least one
from the following equalities is satisfied:}

\begin{equation}
\label{10} (x / u) \cdot (v\backslash yz) = (y(v\backslash x))/ u
\cdot z,
\end{equation}
\begin{equation}
\label{11} (x / u) \cdot (v\backslash (y(v\backslash z)) = (z / u)
\cdot (v\backslash (y(v\backslash x))),
\end{equation}
\begin{equation}
\label{12} (x / u) \cdot (v\backslash yz) = (yz / u) \cdot
(v\backslash x),
\end{equation}
\begin{equation}
\label{13} x(v\backslash ((y / u)z) = ((xz) / u) \cdot
(v\backslash y),
\end{equation}
\begin{equation}
\label{14} x(v\backslash ((y / u) \cdot (v\backslash z)) = (z / u)
\cdot (v\backslash (x(v\backslash y))).
\end{equation}

{\bf Proof:} \emph{Necessity.} Suppose that a quasigroup $(Q,
\cdot)$ is isotopic to an abelian group $(Q, +).$ Let us show that
in a quasigroup the identity (\ref{10}) is satisfied. For this it
is enough to consider a principle isotopy of the form
\begin{equation}
\label{15} x+y=R_{a}^{-1}x\cdot L_{b}^{-1}y.
\end{equation}
From (\ref{15}) follows that
\begin{equation}
\label{16} xy=R_{a}x+L_{b}y.
\end{equation}
Then it is evident that $x/y=R_{a}^{-1}(x+JL_{b}y)$  and
$x\backslash y=L_{b}^{-1}(JR_{a}x+y),$ where $Jx=-x,$ $x\in Q.$

Indeed, denoted $xy=z.$ Then $x=z/y,$ $z=R_{a}(x/y)+L_{b}y,$
$R_{a}(z/ y)=z+JL_{b}y,$ or  $z/y=R_{a}^{-1}(z+JL_{b}y).$

Similarly if $xy=z,$  then  $y=x\backslash z,$
$z=R_{a}x+L_{b}(x\backslash z),$ $L_{b}(x\backslash z)=JR_{a}x+z,$
$x\backslash z=L_{b}^{-1}(JR_{a}x+z).$

Using (\ref{16}) and expressions for enviers operations (/) and
$(\backslash)$, we have:

\begin{align*}
(x/u)\cdot(v\backslash yz)=&R_{a}(x/ u)+L_{b}(v\backslash yz)=\\
=&R_{a}R_{a}^{-1}(x+JL_{b}u)+L_{b}L_{b}^{-1}(JR_{a}v+yz)=\\
=&(x+JL_{b}u)+(JR_{a}v+R_{a}y+L_{b}z)=\\
=& x+JL_{b}u+JR_{a}v+R_{a}y+L_{b}z.
\end{align*}
\begin{align*}
(y(v\backslash x))/ u\cdot z&=R_{a}((y(v\backslash x))/
u)+L_{b}z=\\
&=R_{a}R_{a}^{-1}(y(v\backslash x)+JL_{b}u)+L_{b}z=\\
&=(y \cdot L_{b}^{-1}(JR_{a}v+x)+JL_{b}u)+L_{b}z\\
&=(R_{a}y+(JR_{a}v+x)+JL_{b}u)+L_{b}z=\\
&=R_{a}y+JR_{a}v+x+JL_{b}u+L_{b}z.
\end{align*}

Hence,
\begin{align*}
&(x/ u)(v\backslash yz)=x+JL_{b}u+JR_{a}v+R_{a}y+L_{b}z.\\
&(y(v\backslash x))/ u\cdot z=R_{a}y+JR_{a}v+x+JL_{b}u+L_{b}z.
\end{align*}

As the group $(Q, +)$ is an abelian so the right sides of the
identities  are equal. Therefore in a quasigroup $(Q, \cdot,
\backslash, /)$   the identity (\ref{10}) is satisfied.

\emph{Sufficiency.} Suppose that in a quasigroup $(Q, \cdot, /,
\backslash)$ the identity (\ref{10}) is satisfied . Re-write the
identity (\ref{10}) in the following form:
\begin{align*}
R_{u}^{-1}x\cdot L_{u}^{-1}(yz)=R_{u}^{-1}(y\cdot
L_{u}^{-1}x)\cdot z.
\end{align*}

Put $u=a,$ $v=b$, then
\begin{equation}
\label{17} R_{a}^{-1}x\cdot L_{b}^{-1}(yz)=R_{a}^{-1}(y\cdot
L_{b}^{-1}x)\cdot z.
\end{equation}

The last equality can be written in the form:
\begin{align*}
R_{a}^{-1}x\cdot L_{b}^{-1}(y z)=R_{a}^{-1}(L_{b}^{-1}x \ast
y)\cdot z,
\end{align*}
where  $L_{b}^{-1}x \ast y=y\cdot L_{b}^{-1}x.$

Hence according to Belousov's theorem about four  quasigroups, the
quasigroup $(Q, \ast)$ and also the quasigroup $(Q, \cdot)$ are
isotopic to some group. Therefore, the  principle isotopy $(Q,
+):$ $x+y=R_{a}^{-1}x\cdot L_{b}^{-1}y$ is a group. Using this in
(\ref{17}) we have:
$$
R_{a}^{-1}R_{a}x+L_{b}L_{b}^{-1}(R_{a}y+L_{b}z)=R_{a}R_{a}^{-1}(R_{a}y+L_{b}L_{b}^{-1}x)+L_{b}z\\
$$
or
$$
x+(R_{a}y+L_{b}z)=(R_{a}y+x)+L_{b}z
$$

Replace  $y$ by $R_{a}^{-1}y,$ and $z$ by $L_{b}^{-1}z;$
\begin{align*}
x+(y+z)=(y+x)+z.
\end{align*}

Hence, the  group $(Q, +)$ is  abelian.

The condition of theorem for the identities  (\ref{11}) -
(\ref{14}) are showed analogously.

\endproof

\begin{center}
{\bf III. Chronology development  theory of quasigroups and loops}
\end{center}

Below we show a short chronology of development  of the theory of
quasigroups and loops, also a list  of those who made a
contribution in this theory.

\quad

L. Euler (1783). {\it Latin squares}.\\

R. Moufang (1935). {\it Moufang loops}.\\

A. K. Suschkewitsch (1929,1937). {\it The generalized groups}.\\

B. A. Hausman, O. Ore (1937), G. N. Garrison (1940). {\it Some aspects of the theory of quasigroups}.\\

K. Toyoda, D.S. Murdoch (1939,1941). {\it Medial quasigroups}.\\

A. A. Albert (1943-1944), R. H. Bruck (1944-1946), R. Artzy
(1955-1965), Sh.K. Stein (1957), A. Sade (1960). {\it Basics
algebraic theory quasigoups and loops}.\\

V. D. Belousov (1950-1989). {\it  Foundations of the theory of
quasigroups and loops, Moscow, Nauka, 1967.(in Russian). (First
monograph on structure theory of quasigroups)}.\\

I. Mal'tsev (1970). {\it Connection between Moufang loops and
Mal'tsev's algebras}.\\

T. Evans (1950-1960), M. M. Glukhov,  A. A. Gvaramiya (1969-1971).
{\it The algorithmic problems of the theory of quasigroups and loops}.\\

Y. I. Manin (1972). {\it $CH$-quasigroups}.\\

J. Denes and A. D. Keedwel. (1974). {\it Latin squares and their
applications. Academiai Kiado, Budaest.}\\

J. Smith, O. Chein, H.O. Pflugfelder, T. Kepka, J. Jezek,  P.
Nemec, A. Drapal, G.B. Belyavskaya, V.A. Shcherbacov, K.K.
Shchukin, V.I. Izbash, P.N. Sirbu, F.N. Sokhatsky, E. Kuznetsov,
I.P. Shestakov,  A. Grishkov,  A. Romanowska, W. Dudek, A. Krapez,
G.P. Nagy, M.K. Kinyon, K. Kunen,  J.D. Philips,  P. Vojtechovsky.
{\it The modern a structure theory of quasigroups}.\\

M.A. Akivis, V. V. Goldberg, L.V. Sabinin, P.O. Mikheev. {\it A
quasigroups and a geometry}.

\vspace{7mm}

Presently the theory of quasigroups, as well as other algebraic
structures, is developed in several directions, but among  them,
in our opinion, there are three main, namely:

1) {\it Research the inner nature itself quasigroups;}

2) {\it  A tendency to obtain analogues of known results and
theorems from other algebraic structures;}

3) {\it  Applications of the theory of quasigroups  (cryptology,
differential  geometry, theoretical physics  (Poincare's
quasigroups, Lorentz's quasigroups)}.
\\\\

{\bf Acknowledgments}

I am thankful to professor M.M. Glukhov for  discussion on
receiving  some  results of this article.
\\\\

\begin{center}
{\bf References:}
\end{center}

\begin{enumerate}

\item V.D. Belousov. {\it Foundations of the theory of quasigroups and
loops.} Moscow, Nauka, 1967.(in Russian).

\item  A.Kh. Tabarov.  {\it Linear quasigroups.I }. Preprint,
Central European University, Department Mathematics and its
Applications, Hungary, May 2010, 14 pages.

\item Yu.A. Bakhturin, A.Yu. Olshanskii.  {\it Identical
relations}. Sovremen. problemy matem. Fundament. napravl.,
117-240: in Itogi Nauki i Tekhniki VINITI AN SSSR, 18, Moscow,
1988, (in Russian). Translated in Encyclopaedia of Mathematical
Sciences, Algebra II, Noncommutative Rings. Identities Vol. 18,
Springer, Heidelberg, 1991 Review.

\item G.B. Belyavskaya,  A.Kh. Tabarov.  {\it Characterization of linear and
alinear quasigroups.} Diskretnaya matematika, v.4, no.2, 1992,
p.142-147.(in Russian).

\item G.B. Belyavskaya, A.Kh. Tabarov. {\it
Nuclei and center of linear quasigroups.} Bul. Acad. Stiinte,
Republic Moldova, Matematica, Kishinev (1991), No. 3(6),
p.37-42,(in Russian).

\item  G.B. Belyavskaya, A.Kh. Tabarov.  {\it
One-sided T-quasigroups and irreducible balanced identities.}
Quasigroups and Related Systems. Kishinev, 1994, No. 1, р.8-21.

\item V.D. Belousov.  {\it Balanced identities on quasigroups.} Mat. Sb.,
vol.70(112): no.1, 1966, p.55-97,(in Russian).

\item  F.N. Sokhatski.   {\it Associates and
decompositions of multiary.}  Diss. of Doctor  of Sciences, Kiev,
2006.(in Ukrainian).

\item T. Evans.  {\it The word problem for absract algebras.}
J.London Math.Soc. 1951. 28. No 1, p.64-67.

\item T. Evans.  {\it The isomorphism problem for
some classes of multiplicative systems.} Trans. Amer. Math. Soc.
109(1963), p. 303-312.

\item M.M. Glukhov, A.A.
Gvaramiya.  {\it On  algorithmic problems for some classes of
quasigroups.} DAN USSR, 1967, vol.177, No 1, p.14-16.(in Russian).

\item M.M. Glukhov.  {\it On free products and
algorithmic problems in $R$-varieties of universal algebras.}
Soviet Math. Dokl., vol.11 (1970), No.4,p.957-960.

\item B. Csacany.  {\it On the equivalence of
certain classes of algebraic systems.} Acta Sci. Math. Szeged,
1962, vol.23, p.46-57.

\item J. Jezek  and T. Kepka. {\it Quasigroups, isotopic to a group.}
Commentationes  math. Universitatis Carolinae, 1975, vol.16, No.1,
p.59-76.

\item A.I. Mal'tsev. {\it Identical relations on
varieties of quasigroups.} Mat. Sb., 69(1):p.3-12, 1966,(in
Russian).

\item T. Kepka   and P. Nemec.  {\it $T$-quasigroups. I.} - Acta univ.
Carolin. Math.Phys., 1971, vol.12, No. 1, р.31-39.

\item T. Kepka   and P. Nemec.  {\it $T$-quasigroups. II.} - Acta univ.
Carolin. Math.Phys., 1971, vol. 12, No. 2, р.39-49.

\item  J. Jezek  and T.
Kepka. {\it Varieties of abelian quasigroups.} Czech. Math.
Journal, 1977, vol.27, p.473-503.

\item A.Kh. Tabarov.  {\it A  construction
of free linear quasigroups.} Dokl. Akad. of Sciences of RT, 2005,
vol. XLVIII, No. 11-12, p.22-28.(in Russian).

\item M.M. Glukhov. {\it On free quasigroups of
some varieties and their multiplicative groups.} International
conference devoted to 100-anniversary of A.G.Kurosh.
Abstracts.Moscow, MSU, 2008, p.68.(in Russian).

\item A.Kh. Tabarov.  {\it  Algorithmic solvability of the words problem
for some variety of linear quasigroups.} International Conference.
Mal'tsev  Meeting'09, dedicated to the 100th anniversary of A.I.
Mal'tsev, August 24-28, 2009. Collectin of Abstracts, p.168 (in
Russian).

\item T. Kepka.  {\it Structure of weakly of abelian
quasigroups.} - Czech. Math. Journal, 1978, vol.28, p.181-188.

\item T. Kepka.  {\it Structure of triabelian quasigroups.} - Comment. math.
Univ. Carolinae, 1976, vol.17, p.229-240.

\item P. Nemec.  {\it Commutative
Moufang loops corresponding to linear quasigroups.} - Comment.
math. Univ. Carolinae, 1988, vol. 2, p.303-308.

\item V.A. Shcherbacov.  {\it On linear and inverse quasigroups and their
applications in code theory.} Dissertation of Doctor of Sciences,
247 pages, Chisinau, 2008.

\item K. Toyoda.  {\it On axsioms of linear functions.} - Proc. Imp.
Acad.Tokyo., 1941, vol. 17, p.221-227.

\item V.D. Belousov.  {\it On the structure of distributive
quasigroups.} - Uspekhi Mat. Nauk, 1958, vol.13, No 3, p.235-236
(in Russian).

\item V.D. Belousov,  V.I. Onoi  {\it On loops, isotopic to
leftdistributive quasigroups.} - Math. issled. Kishinev, 1972,
vol. 3(25), с. 135-152.(in Russian).

\item Yu.I. Manin.  {\it Cubic forms.} Moscow,
Nauka , 1972.(in Russian).

\item  G.B. Belyavskaya. {\it Abelian
quasigroups are T-quasigroups.} Quasigroups and Related Systems,
1994, vol.1, no.1, p.1-7.

\item V.D. Belousov. {\it $n$-Ary quasigroups.}
Shtiinsa, Kishinev, 1971.(in Russian).

\item T. Evans.  {\it Abstract mean values.} - Duke math. J., 1963, vol.30,
p.331-347.

\item T. Kepka, M.K. Kinyon  and J.D. Phillips.  {\it The
structure of F-quasigroups.} http://arxiv.org/
abs/math/0510298(2005), 24 pages.

\item V.A. Shcherbacov  {\it On the structure of left and  right F-, SM- and
E-quasigroups.} Journal of Generalized Lie Theory and
Applications, vol.3 (2009), no.3, p.197-259.

\item S.K. Stein.  {\it Left distributive
quasigroups.} - Proc.Amer.Math.Soc., 1959, vol.10, No.4,
p.577-578.

\item V.D. Belousov.  {\it Balanced identities on quasigroups.} Mat. Sb.,
vol.70(112): no.1, 1966, p.55-97,(in Russian).

\item H. Kiechle. {\it Theory of $K$-loops.} Habilitationsschrift. Fachbereich
mathematik der Universitat Hamburg. Hamburg, 1998, (Habilitation
Dissertation).

\item A.A. Gvaramiya.  {\it Aksiomatiziruemie klassi kvazigrupp i mnogosortnaya algebra}.
Dissertation of Doctor of Sciences, Sukhumi, 1985, (in Russian).

\item G.B. Belyavskaya, A.Kh. Tabarov. {\it Identities with permutations
leading to linearity of quasigroups.} Discrete Mathematics and
Applications, 2009, vol. 19, iss. 2, p.173-190.

\item M.I. Kargapolov, Yu.I. Merzlyakov.  {\it Foundations of group theory}. Moscow, Nauka, 1977,
(in Russian).

\item A. Drapal. {\it On multiplication groups of relatively free
quasigroups isotopic to abelian groups.} Czech. Math. J. 55
(130)(2005), p.61-86.

\end{enumerate}

\end{document}